\documentclass[fleqn,preprint,3p,a4paper]{elsarticle}

\usepackage{amssymb,mathrsfs,latexsym,amsfonts}
\usepackage{amsmath}
\usepackage{amsthm}
\usepackage{dcolumn}
\usepackage{endnotes}
\usepackage{tabularx}
\usepackage[matrix,arrow]{xy}
\usepackage{wasysym}

\newtheorem{theorem}{Theorem}[section]
\newtheorem{proposition}[theorem]{Proposition}
\newtheorem{lemma}[theorem]{Lemma}

\theoremstyle{definition}
\newtheorem{definition}[theorem]{Definition}
\newtheorem{example}[theorem]{Example}

\newtheorem{question}[theorem]{Question}

\usepackage{rotating}

\newcommand{\C}{\mathbb{C}}
\newcommand{\R}{\mathbb{R}}

\renewcommand{\Re}{\operatorname{Re}}
\renewcommand{\Im}{\operatorname{Im}}
\newcommand{\gyr}{\operatorname{gyr}}

\journal{Journal}

\begin{document}

\begin{frontmatter}



\title{On two examples of topological gyrogroups\tnoteref{t1}}
\tnotetext[t1]{This research was supported by the National Natural Science Foundation of China (Nos. 12471070), and the Foundation of PhD start-up of Suqian University (Nos. 2024XRC003).}


\author[M. Bao]{Meng Bao}
\ead{mengbao95213@163.com}
\address[M. Bao]{School of Mathematics, Suqian University, Suqian 223800, China}

\author[P. Li]{Piyu Li\corref{mycorrespondingauthor}}
\cortext[mycorrespondingauthor]{Corresponding author.}
\ead{lpy91132006@aliyun.com}
\address[P. Li]{School of Mathematics and Statistics, Xuzhou Institute of Technology, Xuzhou 221018, China}

\begin{abstract}
In this paper, we construct an infinite-dimensional topological gyrogroup that fails to be strongly topological. We equip the countable direct sum of the complex unit disk with a carefully weighted infinite series of component seminorms. By dividing the infinite-dimensional operation into a finite core and a bounded infinite tail, we rigorously prove the joint continuity of the M\"obius addition, confirming the space is a valid topological gyrogroup. Furthermore, we demonstrate that the gyroautomorphisms exhibit unbounded expansion under this real seminorm topology, explicitly preventing the existence of an invariant neighborhood basis at the identity. This provides a clear, analytic counterexample separating general topological gyrogroups from strongly topological ones. Moreover, we construct a topological gyrogroup which is not strongly topological but it contains a non-trivial strong subgyrogroup.
\end{abstract}

\begin{keyword}
Topological gyrogroups, strongly topological gyrogroups, gyr-invariant
\MSC 22A22; 54A20; 20N05; 18A32.

\end{keyword}



%% MSC codes here, in the form:

\end{frontmatter}


\section{Introduction}

The concept of rectifiable spaces was introduced by A.V. Arhangel'ski\v\i\ \cite{Aav80} as non-associative generalizations of topological groups in 1980. Clearly, every topological group is a rectifiable space. It was noted by Uspenski\v\i\ in \cite{UVV,UVV89} that the 7-dimensional sphere $S^{7}$ is an example of a compact rectifiable space, and it does not admit topological group topology. Then Alexandra S. Gul'ko \cite{Alexan} showed that the Hilbert cube $I^{\omega}$ is not rectifiable, while it is a retract of its free topological group. It is interesting to look for some examples of rectifiable spaces, which are not topological groups, except for $S^{7}$.

The concept of a gyrogroup, initially introduced by A. A. Ungar to study the parameterization of the Lorentz group in special relativity, provides a rigorous algebraic framework for non-associative and non-commutative operations \cite{UA2002, UA2005}. Then W. Atiponrat \cite{AW} posed the concept of topological gyrogroups as a generalization of topological groups which is a gyrogroup equipped with a topology such that the binary operation is jointly continuous and the inverse mapping is also continuous. She found that some classical examples of gyrogroups equipped with standard topology are topological gyrogroups but not topological groups, such as the Einstein gyrogroup with the standard topology and the M\"{o}bius gyrogroup with the standard topology. Then, Cai, Lin and He in \cite{CZ} proved that every topological gyrogroup is a rectifiable space, which implies that the classical M\"{o}bius gyrogroup and Einstein gyrogroup equipped with standard topology are rectifiable spaces with infinite cardinal number, which are not topological groups. By further study on the M\"{o}bius gyrogroups, Bao and Lin \cite{BL} introduced the concept of strongly topological gyrogroups requiring the existence of a neighborhood base at the identity that is uniformly invariant under all gyroautomorphisms and found that strongly topological gyrogroups have many properties in common with topological groups. Moreover, they noted that the classical examples of gyrogroups equipped with standard topology are all strongly topological gyrogroups. In order to distinguish the relationship between topological gyrogroups and strongly topological gyrogroups, Bao posed the following question in \cite{BM2024}.

\begin{question}
Is there a topological gyrogroup $G$ but not a strongly topological gyrogroup? Especially with $|G|\geq \omega$.
\end{question}

\begin{question}
Is there a topological gyrogroup $G$ such that $G$ is not a strongly topological gyrogroup and $G$ contains some non-trivial strong subgyrogroups?
\end{question}

In this paper, we construct an explicit infinite-dimensional counterexample, answering this question in the negative.

\section{Preliminary}

Throughout this paper, all topological spaces are assumed to be Hausdorff, unless otherwise is explicitly stated. Let $\mathbb{N}$ be the set of all positive integers and $\omega$ the first infinite ordinal. The readers may consult \cite{AA2008, E, linbook1, UA} for notation and terminology not explicitly given here. Next we recall some definitions and facts.

\begin{definition}(\cite{AW})
Let $G$ be a nonempty set, and let $\oplus: G\times G\rightarrow G$ be a binary operation on $G$. Then the pair $(G, \oplus)$ is called a {\it magma or groupoid}. A function $f$ from a groupoid $(G_{1}, \oplus_{1})$ to a groupoid $(G_{2}, \oplus_{2})$ is called a {\it groupoid homomorphism} if $f(x\oplus_{1}y)=f(x)\oplus_{2} f(y)$ for any elements $x, y\in G_{1}$. Furthermore, a bijective groupoid homomorphism from a groupoid $(G, \oplus)$ to itself will be called a {\it groupoid automorphism}. We write $\mbox{Aut}(G, \oplus)$ for the set of all automorphisms of a groupoid $(G, \oplus)$.
\end{definition}

\begin{definition}(\cite{UA})
Let $(G, \oplus)$ be a groupoid. The system $(G,\oplus)$ is called a {\it gyrogroup}, if its binary operation satisfies the following conditions:

\smallskip
$(G1)$ There exists a unique identity element $0\in G$ such that $0\oplus a=a=a\oplus0$ for all $a\in G$.

\smallskip
$(G2)$ For each $x\in G$, there exists a unique inverse element $\ominus x\in G$ such that $\ominus x \oplus x=0=x\oplus (\ominus x)$.

\smallskip
$(G3)$ For all $x, y\in G$, there exists $\gyr[x, y]\in \mbox{Aut}(G, \oplus)$ with the property that $x\oplus (y\oplus z)=(x\oplus y)\oplus \gyr[x, y](z)$ for all $z\in G$.

\smallskip
$(G4)$ For any $x, y\in G$, $\gyr[x\oplus y, y]=\gyr[x, y]$.
\end{definition}

Clearly, every group is a gyrogroup, where the groupoid automorphism is the identity function.

\begin{definition}(\cite{ST})
Let $(G,\oplus)$ be a gyrogroup. A nonempty subset $H$ of $G$ is called a {\it subgyrogroup}, denoted
by $H\leq G$, if $H$ forms a gyrogroup under the operation inherited from $G$ and the restriction of $\gyr[a,b]$ to $H$ is an automorphism of $H$ for all $a,b\in H$.

\smallskip
Furthermore, a subgyrogroup $H$ of $G$ is said to be an {\it $L$-subgyrogroup}, denoted
by $H\leq_{L} G$, if $\gyr[a, h](H)=H$ for all $a\in G$ and $h\in H$.
\end{definition}

\begin{definition}(\cite{BX2022})
A subgyrogroup $H$ of a topological gyrogroup $G$ is called {\it strong subgyrogroup} if for any $x,y\in G$, we have $\gyr[x,y](H)=H$.
\end{definition}

Obviously, every strong subgyrogroup is an $L$-subgyrogroup.

\begin{lemma}\emph{(\cite{UA})}\label{a}.
Let $(G, \oplus)$ be a gyrogroup. Then for any $x, y, z\in G$, we obtain the following:

\smallskip
1. $(\ominus x)\oplus (x\oplus y)=y$;

\smallskip
2. $(x\oplus (\ominus y))\oplus \gyr[x, \ominus y](y)=x$;

\smallskip
3. $(x\oplus \gyr[x, y](\ominus y))\oplus y=x$;

\smallskip
4. $\gyr[x, y](z)=\ominus (x\oplus y)\oplus (x\oplus (y\oplus z))$;

\smallskip
5. $(\ominus x\oplus y)\oplus \gyr[\ominus x, y](\ominus y\oplus z)=\ominus x\oplus z$;

\smallskip
6. $(x\oplus y)\oplus z=x\oplus (y\oplus \gyr[y,x](z))$;

\smallskip
7. $A\gyr[x,y]=\gyr[Ax, Ay]A$ \mbox{ for any automorphism} $A$ \mbox{ on } $G$;

\smallskip
8. $\gyr[\ominus x, \ominus y]=\gyr[x,y]$;

\smallskip
9. $\gyr^{-1}[x, y]=\gyr[y, x]$;

\smallskip
10. $\ominus (x\oplus y)=\gyr[x,y](\ominus y\ominus x)$.
\end{lemma}

\begin{proposition}\label{subset=}\emph{(\cite{ST})}
Let $G$ be a gyrogroup and let $X\subseteq G$. Then the following are equivalent:

\emph{(1)} $\gyr[x,y](X)\subseteq X$ for all $x,y\in G$;

\emph{(2)} $\gyr[x,y](X)=X$ for all $x,y\in G$;
\end{proposition}

\begin{definition}(\cite{AW})
A triple $(G, \tau, \oplus)$ is called a {\it topological gyrogroup} if the following statements hold:

\smallskip
(1) $(G, \tau)$ is a topological space.

\smallskip
(2) $(G, \oplus)$ is a gyrogroup.

\smallskip
(3) The binary operation $\oplus: G\times G\rightarrow G$ is jointly continuous while $G\times G$ is endowed with the product topology, and the operation of taking the inverse $\ominus (\cdot): G\rightarrow G$, i.e. $x\rightarrow \ominus x$, is also continuous.
\end{definition}

It is clear that every topological group is a topological gyrogroup. The following example is a topological gyrogroup, but not a topological group.

\begin{example}(\cite[Example 3]{AW})
The Einstein gyrogroup with the standard topology is a topological gyrogroup but not a topological group.
\end{example}

By further research on topological gyrogroups, Bao and Lin in \cite{BL} introduced the concept of strongly topological gyrogroups.

\begin{definition}{\rm (\cite{BL})}\label{d11}
Let $G$ be a topological gyrogroup. We say that $G$ is a {\it strongly topological gyrogroup} if there exists a neighborhood base $\mathscr U$ of $0$ such that, for every $U\in \mathscr U$, $\gyr[x, y](U)=U$ for any $x, y\in G$. For convenience, we say that $G$ is a strongly topological gyrogroup with neighborhood base $\mathscr U$ of $0$.
\end{definition}

The following example is a strongly topological gyrogroup and it is not a topological group.

\begin{example}(\cite[Example 3.1]{BL})
The M\"{o}bius gyrogroup with the standard topology is a strongly topological gyrogroup but not a topological group.
\end{example}

Then, based on the M\"{o}bius gyrogroup, we construct a topological gyrogroup such that it does not have any gyr-invariant base at the identity element, that is, it is not a strongly topological gyrogroup.

\section{An infinite-dimensional topological gyrogroup but not strongly topological}

\subsection{Construction of the Space}
Let $D = \{ z \in \mathbb{C} : |z| < 1 \}$ be the open complex unit disk equipped with the standard M\"obius addition:
$$x \oplus_M y = \frac{x + y}{1 + \bar{x}y}.$$
The space $(D, \oplus_M)$ is a well-known gyrogroup. The gyration operator in $D$ is given by:
$$\gyr_M[x,y](z) = \frac{x\bar{y} + 1}{\bar{x}y + 1} z.$$
Because the coefficient $\frac{x\bar{y} + 1}{\bar{x}y + 1}$ has a modulus of $1$, $\gyr_M[x,y]$ acts as a rigid rotation on $D$.

We define our space $G$ as the finite-support infinite direct sum of $D$:
$$G = \bigoplus_{n=1}^\infty D_n.$$
An element $z \in G$ is a sequence $(z_1, z_2, z_3, \dots)$ where $z_n \in D_n=D$, and only finitely many $z_n$ are non-zero. Operations are defined coordinate-wise. Thus, $(G, \oplus)$ strictly satisfies all gyrogroup axioms.

\begin{definition}[The Weighted Component Seminorm]
For each integer $n \ge 1$, we define a mapping $p_n: \C \to \R_{\ge 0}$ by breaking a complex number into its real and imaginary parts with an asymmetric weight:
$p_n(w) = |\Re(w)| + n|\Im(w)|.$
\end{definition}

\begin{proposition} \label{prop:seminorm}
The mapping $p_n(w)$ is a well-defined real seminorm on the complex vector space $\C$ (treated as a real vector space).
\end{proposition}
\begin{proof}
Let $w = x + iy$ and $w' = x' + iy'$, where $x, y, x', y'\in \R$. The non-negativity $p_n(w) \ge 0$ is trivial since $n \ge 1$. For subadditivity (the triangle inequality), we have:
\begin{align*}
    p_n(w + w') &= |x + x'| + n|y + y'| \\
    &\le (|x| + |x'|) + n(|y| + |y'|) = p_n(w) + p_n(w').
\end{align*}
For absolute homogeneity with respect to any \textit{real} scalar $c \in \R$:
$$p_n(cw) = |cx| + n|cy| = |c|(|x| + n|y|) = |c|p_n(w).$$
Thus, $p_n$ satisfies all axioms of a real seminorm. Note that it does not satisfy complex absolute homogeneity (e.g., $p_n(iw) \neq |i|p_n(w)$), which plays a critical role in our main result.
\end{proof}

\begin{definition}[The Infinite-Dimensional Metric Space and Gyrogroup]
Let $V = \bigoplus_{n=1}^\infty \C_n$ (where $\C_n=\C$) be the complex vector space of all complex sequences with only finitely many non-zero terms. We define a functional $P(\cdot)$ on $V$ by the infinite series:
$$
    P(z) = \sum_{n=1}^\infty p_n(z_n), \quad \forall z \in V.
$$
By Proposition \ref{prop:seminorm}, each individual $p_n$ is a real norm on $\C$. Since every sequence in $V$ has only finitely many non-zero terms, the series $P(z)$ always reduces to a finite sum. A finite sum of real norms trivially satisfies positive definiteness, subadditivity, and real absolute homogeneity. Consequently, $P(z)$ constitutes a valid real norm on the ambient vector space $V$.

This norm naturally induces a translation-invariant metric on $V$ defined by standard vector subtraction:
$$
    d(x, y) = P(x - y).
$$

Now, let $D = \{ z \in \C : |z| < 1 \}$ be the open complex unit disk, and define the subset $G = \bigoplus_{n=1}^\infty D_n \subset V$, where $D_n=D$. The set $G$ forms a strictly defined gyrogroup under component-wise M\"obius addition. We endow the gyrogroup $G$ with the relative metric topology inherited from $(V, d)$. Under this strictly metrizable topology, the basic open neighborhoods of any $x \in G$ are $B_\epsilon(x) = \{z \in G : P(z - x) < \epsilon\}$. Sequence convergence $z^{(m)} \to x$ in $G$ is precisely equivalent to $P(z^{(m)} - x) \to 0$ as $m \to \infty$.
\end{definition}

% ==========================================
\subsection{Continuity of the M\"obius Addition}
% ==========================================

To prove that $(G, \tau_P, \oplus)$ is a topological gyrogroup, we must establish the joint continuity of $\oplus$.

\begin{lemma}[Topological Equivalence] \label{lem:equivalence}
For any fixed integer $n \ge 1$, $p_n(w)$ is topologically equivalent to the Euclidean modulus $|w|$, satisfying:
$|w| \le p_n(w) \le (n+1)|w|.$
\end{lemma}
\begin{proof}
Let $w = u + iv$. For the lower bound, $$p_n(w) = |u| + n|v| \ge |u| + |v| \ge \sqrt{u^2 + v^2} = |w|.$$ For the upper bound, $p_n(w) \le |w| + n|w| = (n+1)|w|$.
\end{proof}

\begin{theorem}
The space $(G, \tau_P, \oplus)$ is a topological gyrogroup.
\end{theorem}
\begin{proof}
Let $x, y \in G$. Since $x$ and $y$ are elements of the direct sum space, there exists a finite integer $K \ge 1$ such that $x_n = y_n = 0$ for all $n > K$.

Consider sequences $u^{(m)} \to x$ and $v^{(m)} \to y$ as $m \to \infty$ under the norm $P(\cdot)$. We need to prove that $P(u^{(m)} \oplus v^{(m)} - x \oplus y) \to 0$. We split the infinite series into two parts:
\begin{equation*}
    P(u^{(m)} \oplus v^{(m)} - x \oplus y) = \sum_{n=1}^K p_n(u^{(m)}_n \oplus v^{(m)}_n - x_n \oplus y_n) + \sum_{n=K+1}^\infty p_n(u^{(m)}_n \oplus v^{(m)}_n).
\end{equation*}

\textbf{Part 1: The Finite Core ($n \le K$).} \\
For any fixed index $n \in \{1, \dots, K\}$, since all terms in the sum $P(\cdot)$ are non-negative, we have the elementary bound:
$$
    p_n(u^{(m)}_n - x_n) \le \sum_{k=1}^\infty p_k(u^{(m)}_k - x_k) = P(u^{(m)} - x).
$$
Because $P(u^{(m)} - x) \to 0$ as $m \to \infty$, it follows strictly that $p_n(u^{(m)}_n - x_n) \to 0$. By Lemma \ref{lem:equivalence}, we deduce the convergence in the standard Euclidean metric:
$$
    |u^{(m)}_n - x_n| \le p_n(u^{(m)}_n - x_n) \to 0.
$$
Thus, $u^{(m)}_n \to x_n$ in $\C$. By an identical argument, $v^{(m)}_n \to y_n$ in $\C$.

Furthermore, since $x, y \in G = \bigoplus_{n=1}^\infty D$, their components natively reside in the open unit disk, meaning $|x_n| < 1$ and $|y_n| < 1$. Let us define the complex M\"obius addition as a mapping $f: D \times D \to D$, where $f(u, v) = \frac{u+v}{1+\bar{u}v}$. At the specific point $(x_n, y_n)$, the denominator is strictly bounded away from zero because:
$$
    |1 + \bar{x}_n y_n| \ge 1 - |x_n||y_n| > 0.
$$
Being a rational function with a non-vanishing denominator on $D \times D$, the mapping $f$ is jointly continuous with respect to the standard Euclidean metric. Recall from Lemma \ref{lem:equivalence} that convergence in the seminorm $p_n$ is strictly equivalent to convergence in the Euclidean metric. Thus, our established component-wise convergence implies the Euclidean convergence $u^{(m)}_n \to x_n$ and $v^{(m)}_n \to y_n$. By invoking the sequential continuity of $f$, the limit operator smoothly commutes with the function evaluation, yielding:
\begin{align}
    \lim_{m \to \infty} |u^{(m)}_n \oplus v^{(m)}_n - x_n \oplus y_n|
    &= \left| \lim_{m \to \infty} f(u^{(m)}_n, v^{(m)}_n) - f(x_n, y_n) \right| \nonumber \\
    &= \left| f\left( \lim_{m \to \infty} u^{(m)}_n, \lim_{m \to \infty} v^{(m)}_n \right) - f(x_n, y_n) \right| \nonumber \\
    &= |f(x_n, y_n) - f(x_n, y_n)| = 0.
\end{align}

To translate this Euclidean convergence back to our weighted norm $p_n$, we apply the upper bound from Lemma \ref{lem:equivalence}:
$$
    p_n(u^{(m)}_n \oplus v^{(m)}_n - x_n \oplus y_n) \le (n+1) |u^{(m)}_n \oplus v^{(m)}_n - x_n \oplus y_n| \to 0.
$$
Since this convergence holds for each individual $n \in \{1, \dots, K\}$, and $K$ is a fixed finite integer, the finite sum of these limits converges to zero:
$$
    \lim_{m \to \infty} \sum_{n=1}^K p_n(u^{(m)}_n \oplus v^{(m)}_n - x_n \oplus y_n) = \sum_{n=1}^K 0 = 0.
$$

\textbf{Part 2: The Infinite Tail ($n > K$).} \\
We must ensure the asymmetric weights $n$ in $p_n(w)$ do not cause divergence in the tail. By our selection of $K$, $x_n = 0$ for all $n > K$. Consequently, for any $n > K$, we can establish a uniform upper bound:
\begin{equation*}
    |u^{(m)}_n| \le p_n(u^{(m)}_n) \le \sum_{j=K+1}^\infty p_j(u^{(m)}_j) \le P(u^{(m)} - x).
\end{equation*}
Since $P(u^{(m)} - x) \to 0$ and $P(v^{(m)} - y) \to 0$ as $m \to \infty$, there exists a sufficiently large integer $M$ such that for all $m > M$, we have $P(u^{(m)} - x) < \frac{1}{2}$ and $P(v^{(m)} - y) < \frac{1}{2}$. This strictly guarantees that $|u^{(m)}_n| < \frac{1}{2}$ and $|v^{(m)}_n| < \frac{1}{2}$ simultaneously for all $n > K$.

For the remainder of the proof, we assume $m > M$. Let $$u^{(m)}_n = x_1 + i y_1 \quad \text{and } v^{(m)}_n = x_2 + i y_2.$$

The M\"obius addition yields $u^{(m)}_n \oplus v^{(m)}_n = (\alpha + i\beta)(X + iY)$, where $X = x_1+x_2$, $Y = y_1+y_2$, and the complex scalar is $\frac{1}{1+\bar{u}^{(m)}_n v^{(m)}_n} = \alpha + i\beta$.
The denominator modulus squared is $D = |1+\bar{u}^{(m)}_n v^{(m)}_n|^2 \ge (1 - 1/4)^2 = 9/16$.

Applying the triangle inequality to the real and imaginary components yields:
\begin{align} \label{eq:tail_expand}
    p_n(u^{(m)}_n \oplus v^{(m)}_n) &= |\alpha X - \beta Y| + n|\beta X + \alpha Y| \nonumber \\
    &\le |\alpha|(|X| + n|Y|) + |\beta||Y| + n|\beta||X| \nonumber \\
    &= |\alpha| p_n(u^{(m)}_n + v^{(m)}_n) + |\beta||Y| + n|\beta||X|.
\end{align}

Crucially, the imaginary part of the denominator is $\beta = \frac{-(x_1 y_2 - x_2 y_1)}{D}$. Multiplying by $n$ allows us to pair the weight precisely with the $y$-components:
\begin{equation*}
    n|\beta| \le \frac{1}{D} \big( |x_1| (n|y_2|) + |x_2| (n|y_1|) \big).
\end{equation*}
Since $|x_1|, |x_2| < \frac{1}{2}$, and observing $n|y_1| \le p_n(u^{(m)}_n)$ and $n|y_2| \le p_n(v^{(m)}_n)$, we obtain $$n|\beta| \le \frac{1}{2D} \big( p_n(u^{(m)}_n) + p_n(v^{(m)}_n) \big).$$
Given $|X| < 1$, the cross-term is safely bounded: $n|\beta||X| \le \frac{1}{2D} \big( p_n(u^{(m)}_n) + p_n(v^{(m)}_n) \big)$.

To rigorously bound the non-linear cross terms generated by the M\"obius addition, we analyze the quantity $|\beta||Y|$, where $\beta= \frac{y_1 x_2 - x_1 y_2}{D}$
and $Y=y_1 + y_2$.

Since $m > M$, we possess the uniform coordinate bounds $|x_1| < \frac{1}{2}$ and $|x_2| < \frac{1}{2}$. Applying the triangle inequality yields:
\begin{equation*}
    |\beta| = \frac{|y_1 x_2 - x_1 y_2|}{D} \le \frac{|y_1||x_2| + |x_1||y_2|}{D} < \frac{1}{2D}(|y_1| + |y_2|).
\end{equation*}
Crucially, the uniform bounds $|y_1| < \frac{1}{2}$ and $|y_2| < \frac{1}{2}$ strictly guarantee that $(|Y|\leq |y_1| + |y_2|) < 1$. It follows that
\begin{equation*}
    |\beta||Y| < \frac{1}{2D} (|y_1| + |y_2|).
\end{equation*}
Recalling the definition of our sequence seminorm, $p_n(w) = |\Re(w)| + n|\Im(w)| \ge |\Im(w)|$ for all $n \ge 1$. Consequently, $|y_1| \le p_n(u^{(m)}_n)$ and $|y_2| \le p_n(v^{(m)}_n)$. Then
\begin{equation*} \label{eq:beta_Y_bound}
    |\beta||Y| < \frac{1}{2D} \Big( p_n(u^{(m)}_n) + p_n(v^{(m)}_n) \Big).
\end{equation*}

The real part $\alpha = \frac{1 + x_1 x_2 + y_1 y_2}{D}$ satisfies $|\alpha| \le \frac{1 + |u||v|}{D} < \frac{5/4}{D}$. Summing all bounded components in:
\begin{equation*}
    p_n(u^{(m)}_n \oplus v^{(m)}_n) \le \left( \frac{5/4}{D} + \frac{1}{2D} + \frac{1}{2D} \right) \big( p_n(u^{(m)}_n) + p_n(v^{(m)}_n) \big) = \frac{9/4}{D} \big( p_n(u^{(m)}_n) + p_n(v^{(m)}_n) \big).
\end{equation*}
Given $D \ge 9/16$, the constant $\frac{9/4}{D} \le 4$. Therefore, $$p_n(u^{(m)}_n \oplus v^{(m)}_n) \le 4 \big( p_n(u^{(m)}_n) + p_n(v^{(m)}_n) \big).$$

Summing over $n > K$, the tail is bounded by $4 P(u^{(m)} - x) + 4 P(v^{(m)} - y)$, which strictly converges to $0$. Combining both parts, $P(u^{(m)} \oplus v^{(m)} - x \oplus y) \to 0$.
\end{proof}

% ==========================================
\subsection{Failure of the Strongly Topological Property}
% ==========================================

\begin{theorem}
The topological gyrogroup $(G, \tau_P, \oplus)$ is not strongly topological.
\end{theorem}
\begin{proof}
Assume towards contradiction that $G$ is strongly topological. Let $V_\epsilon(0) = \{ z \in G : P(z) < \epsilon \}$. There must exist a uniformly invariant neighborhood base at $0$. Thus, there exists some base neighborhood $U$ and an $\epsilon > 0$ such that $V_\epsilon(0) \subseteq U \subseteq V_1(0)$, and for all $a, b \in G$, $\gyr[a,b](U) = U$. This strictly implies that $\gyr[a,b](V_\epsilon(0)) \subseteq V_1(0)$.

We now construct a counterexample. Choose a coordinate dimension $M$ sufficiently large such that $M > \frac{4}{\epsilon}$. Let $z \in G$ be defined such that $z_M = \frac{\epsilon}{2}$ (purely real), and $z_n = 0$ for all $n \neq M$.
Clearly, $P(z) = p_M(z_M) = \left| \frac{\epsilon}{2} \right| + M \cdot 0 = \frac{\epsilon}{2} < \epsilon$. Thus, $z \in V_\epsilon(0)$.

We find elements $a, b \in G$ to induce a $90^\circ$ phase rotation exclusively at coordinate $M$. We require the gyration phase factor to be $i$:
$$
    \frac{a_M \bar{b}_M + 1}{\bar{a}_M b_M + 1} = i \implies \text{Arg}(a_M \bar{b}_M + 1) = \frac{\pi}{4}.
$$
Let $a_M \bar{b}_M = -0.5 + 0.5i$. We systematically choose $b_M = (0.5)^{1/4} \approx 0.84 < 1$ (real) and $a_M = (-0.5 + 0.5i) / (0.5)^{1/4}$. The modulus $|a_M|^2 = \frac{0.5}{\sqrt{0.5}} = \sqrt{0.5} < 1$. We set all other coordinates of $a$ and $b$ to $0$. Thus, $a, b \in G$ are well-defined elements within the unit disks.

Let $w = \gyr[a,b]z$. The operator executes a component-wise operation, rotating $z_M$ by exactly $90^\circ$:
$$
    w_M = i \cdot z_M = i \frac{\epsilon}{2}.
$$
For all $n \neq M$, $a_n=b_n=0$, so $\gyr[0,0]z_n = z_n = 0$. We calculate the topological size of $w$:
$$
    P(w) = p_M(w_M) = |\Re(w_M)| + M \cdot |\Im(w_M)| = 0 + M \cdot \frac{\epsilon}{2}.
$$
Since we explicitly chose $M > \frac{4}{\epsilon}$, we obtain:
$$
    P(w) > \frac{4}{\epsilon} \cdot \frac{\epsilon}{2} = 2 > 1.
$$
Therefore, $w \notin V_1(0)$.

This demonstrates that a valid gyration maps an element $z \in V_\epsilon(0)$ strictly outside of $V_1(0)$, breaking the required condition $\gyr[a,b](V_\epsilon(0)) \subseteq V_1(0)$. This unbounded rotational stretch explicitly breaks uniform equicontinuity, proving that no gyration-invariant neighborhood base exists. Consequently, $(G, \tau_P, \oplus)$ is not a strongly topological gyrogroup.
\end{proof}

\section{A topological gyrogroup contains a non-trivial strong subgyrogroup}

The following result shows that even if a topological gyrogroup is metrizable, locally compact and connected space, it may not be a strongly topological gyrogroup.

\begin{proposition}
There exists a metrizable, locally compact, connected Hausdorff topological gyrogroup such that there is no neighborhood base at the identity consisting of $\gyr$-invariant neighborhoods and contains a non-trivial strong subgyrogroup.
\end{proposition}

\begin{proof}
Let $G=\mathbb R^3$, equipped with the usual Euclidean topology. For $x=(a,b,c),y=(p,q,r),z=(u,v,w)\in G$, define the gyrooperation
$$x\oplus y=\big(a+p,\;b+q,\;c+r+(a+2p)(aq-bp)\big),$$
the inverse $\ominus x=(-a,-b,-c)$, and the identity element $e=(0,0,0)$. The gyrator is defined by
$$
\gyr[x,y](u,v,w)=\big(u,\;v,\;w+3(aq-bp)u\big)=S_{3(aq-bp)}(u,v,w),
$$
where $S_t(u,v,w)=(u,v,w+tu)$ denotes the shear transformation.

1. \textbf{Identity axiom.}
$$
e\oplus x = (0+a,0+b,0+c+(0+2a)(0\cdot b - a\cdot 0))=(a,b,c)=x,
$$
$$
x\oplus e = (a+0,b+0,c+0+(a+0)(a\cdot 0-b\cdot 0))=(a,b,c)=x.
$$

2. \textbf{Inverse axiom.}
$$
(\ominus x)\oplus x = (-a+a,-b+b,-c+c+(-a+2a)((-a)b - (-b)a))=(0,0,0)=e,
$$
$$
x\oplus (\ominus x) = (a-a,b-b,c-c+(a-2a)(a(-b)-b(-a)))=(0,0,0)=e.
$$

3. \textbf{Left gyroassociative law.}

Set $F(a,b,p,q)=(a+2p)(aq-bp)$. We compute
$$
y\oplus z = (p+u,q+v,r+w+F(p,q,u,v)),
$$
$$
x\oplus(y\oplus z)=\big(a+p+u,\;b+q+v,\;c+r+w + F(p,q,u,v)+F(a,b,p+u,q+v)\big).
$$
On the other hand,
$$
x\oplus y=(a+p,b+q,c+r+F(a,b,p,q)),
\mbox{and} \gyr[x,y]z=(u,v,w+3(aq-bp)u),
$$
$$
(x\oplus y)\oplus \gyr[x,y]z
=\big(a+p+u,\;b+q+v,\; c+r+F(a,b,p,q)+w+3(aq-bp)u + F(a+p,b+q,u,v)\big).
$$
The first two coordinates coincide trivially.

Since $$a(q+v)-b(p+u)=(aq-bp)+(av-bu)$$ and $$(a+p)v-(b+q)u=(av-bu)+(pv-qu),$$ we obtain
$$
\begin{aligned}
&F(p,q,u,v)+F(a,b,p+u,q+v) -F(a,b,p,q)-F(a+p,b+q,u,v)\\
&=(p+2u)(pv-qu)  +(a+2p+2u)\bigl[(aq-bp)+(av-bu)\bigr]\\
&\quad -(a+2p)(aq-bp)  -(a+p+2u)\bigl[(av-bu)+(pv-qu)\bigr]\\
&=-a(pv-qu)+2u(aq-bp)+p(av-bu)\\
&=-apv+aqu+2aqu-2bpu+apv-bpu\\
&=3(aq-bp)u.
\end{aligned}
$$

Therefore,
$$
F(p,q,u,v)+F(a,b,p+u,q+v) = F(a,b,p,q)+F(a+p,b+q,u,v)+3(aq-bp)u,
$$
which deduces that
 $$
x\oplus(y\oplus z) = (x\oplus y)\oplus \gyr[x,y] z.
$$Hence, the left gyroassociative law holds.

4. \textbf{Each gyrator is an automorphism.}
We verify that $S_t$ preserves the gyrooperation $\oplus$:
$$
S_t(x\oplus y) =\big(a+p,b+q,c+r+F(a,b,p,q)+t(a+p)\big),
$$
$$
S_t x \oplus S_t y =(a+p,b+q,(c+ta)+(r+tp)+F(a,b,p,q))
=\big(a+p,b+q,c+r+F(a,b,p,q)+t(a+p)\big).
$$
Thus $S_t(x\oplus y)=S_t x \oplus S_t y$. Since $S_t$ is bijective and linear, it is an automorphism of $(G,\oplus)$. As $\gyr[x,y]=S_{3(aq-bp)}$, we have $\gyr[x,y]\in\operatorname{Aut}(G,\oplus)$ for all $x,y\in G$.

5. \textbf{Left loop property (gyrator identity):} $\gyr[x\oplus y,y]=\gyr[x,y]$.
$$
\gyr[x\oplus y,y] = S_{3\left((a+p)q-(b+q)p\right)} =S_{3(aq-bp)}=\gyr[x,y],
$$
where we use the identity $(a+p)q-(b+q)p = aq-bp$.

The above axioms confirm that $(G,\oplus)$ is a gyrogroup algebraically.

6. \textbf{Topological gyrogroup conditions.}
Equipped with the Euclidean topology on $\mathbb R^3$:
The gyrooperation $\oplus:G\times G\to G$ is continuous, since all its coordinate functions are polynomials;
The inversion map $\ominus: x\mapsto -x$ is linear and continuous;
The gyrator map $\gyr[\cdot,\cdot](\cdot):G\times G\times G\to G$ is continuous, as every component is a polynomial.

Therefore, $(G,\oplus,\tau_{\text{euclid}})$ is a Hausdorff, metrizable, locally compact, connected topological gyrogroup.

7. \textbf{There is no neighborhood base at $e$ consisting of $\gyr$-invariant neighborhoods.}
A neighborhood $W$ of $e$ is called $\gyr$-invariant provided $\gyr[x,y](W)=W$ for every $x,y\in G$.
Take $x=(1,0,0),\;y=(0,1,0)$. Then $aq-bp = 1\cdot 1 - 0\cdot 0 = 1$, so
$$
\gyr[x,y]=S_{3},\qquad T:=S_3(u,v,w)=(u,v,w+3u).
$$
Suppose for contradiction that there exists a neighborhood base at $e$ formed by $\gyr$-invariant neighborhoods. Consider the open neighborhood $U=\{(u,v,w):|w|<1\}$ of $e$. By the definition of neighborhood base, there exists a $\gyr$-invariant neighborhood $V\subseteq U$. Since $V$ is a neighborhood of the identity element, there exists $\varepsilon>0$ such that $p_0=(\varepsilon,0,0)\in V$. By invariance under $T$,
$$
T^n(\varepsilon,0,0)=(\varepsilon,0, 3n\varepsilon)\in V,\quad \forall n\in\mathbb N.
$$
Choose $n\in\mathbb N$ large enough so that $3n\varepsilon>1$. Then $T^n p_0\notin U$, while $T^n p_0\in V\subseteq U$, a contradiction.

Hence there is no $\gyr$-invariant neighborhood base at the identity.
$\square$

8. \textbf{There exists a non-trivial strong subgyrogroup $H$ in $G$.}

Let $H=\{(0,b,c):b,c\in\mathbb R\}$. For any $h_1=(0,b,c)$ and $h_2=(0,q,r)$ in $H$, the correction term in the third coordinate of $h_1\oplus h_2$ vanishes, since $(0+2\cdot 0)(0\cdot q-b\cdot 0)=0$. Thus $h_1\oplus h_2=(0,b+q,c+r)\in H$. Moreover, $e=(0,0,0)\in H$ and $\ominus h_1=(0,-b,-c)\in H$. Hence $H$ is a subgyrogroup of $G$, and its inherited operation is coordinatewise addition.

To verify that $H$ is strong, let $x=(a,b,c)$ and $y=(p,q,r)$ be arbitrary elements of $G$, and let $h=(0,v,w)\in H$. Then $\gyr[x,y](h)=(0,v,w+3(aq-bp)\cdot 0)=(0,v,w)=h$. Thus every gyration of $G$ fixes $H$ pointwise. In particular, $\gyr[x,y](H)=H$ for all $x,y\in G$.

Since $(0,1,0)\in H\setminus\{e\}$ and $(1,0,0)\notin H$, the subgyrogroup $H$ is nontrivial and proper. It is also closed, being the inverse image of $\{0\}$ under the first coordinate projection.

\end{proof}

\end{document}